\documentclass[12pt]{amsart}
\title{A little Help from my Friends}
\author{%
Anthony G. O'Farrell}
\address{%
Matematics Department\\
NUI, Maynooth\\
Co. Kildare\\
Ireland}
\date{\today}

\usepackage{amsmath,latexsym,epsfig,graphics,amsfonts,amssymb}

\newtheorem{theorem}{Theorem}[section]

\newtheorem{corollary}[theorem]{Corollary}

\begin{document}
\maketitle

\def\scr#1{{\mathcal#1}}

This article is based on a talk given at a
one-day meeting in NUI, Maynooth on the Fourth of April, 2008,
held to honour David Walsh and Richard Watson.

\section{Introduction} This is a tribute to my dear colleagues and friends
David Walsh and Richard Watson, who were here before me in Maynooth, and who laboured with me in the day and the heat. They cheerfully shouldered with me a teaching load that would, apparently, kill the academics of today. The teaching load required, in order to ensure that our students were adequately trained, continued to be a problem until the presidency of M\'ichael Ledwith, in the early
1990's. It was not easy to pursue research while giving 275 lectures a year, but they gave it their best.  At distinct times, both helped me in my investigations. David was sound on complex analysis and hard analysis, so we joined forces to tackle some problems that required technical estimates for integral kernels that solve the $\bar\partial$-problem. Richard had a sound background in
algebra, and he worked with me on problems that could be addressed using algebras of smooth functions.
Richard became my Ph.D. student, after a while, and then, after graduating, continued to work with me for a few years. My period of active collaboration with 
David was in the early eighties, and with Richard the nineties. Recently, both have taken some interest in my reversibility project.

Most of the sources referred to in what follows will be 
found in the references cited in our joint papers, which are listed 
in the bibliography below.

\section{A Way to think of Complex Analysis}
Holomorphic functions are the solutions to the $\bar\partial$ equation
$$ \bar\partial f=0,$$
where
$$ \bar\partial f =\frac{\partial f}{\partial\bar z}\,d\bar z,\qquad
\frac{\partial f}{\partial\bar z} = \frac12\left\{
\frac{\partial f}{\partial x}
+i\frac{\partial f}{\partial y}
\right\}.
$$
We shall refer to both $\bar\partial$ and $\frac{\partial}{\partial\bar z}$
as \lq\lq the $\bar\partial$ operator" (pronounced {\it $d$-bar operator}), 
as convenient.

The $\bar\partial$ operator is skew:
$$\int_{\mathbb{C}}\phi\frac{\partial\psi}{\partial\bar z}\,dx\,dy = 
-\int_{\mathbb{C}}\psi\frac{\partial\phi}{\partial\bar z}\,dx\,dy ,$$
whenever $\phi$ and $\psi$ belong to the space $\scr D$ of $C^\infty$
complex-valued functions on $\mathbb{C}$, having compact support.

The adjoint operator acts on distributions:
$$\left\langle\phi,\left(\frac{\partial}{\partial\bar z}\right)^*f\right\rangle  
=
\left\langle\frac{\partial}{\partial\bar z}\phi,f\right\rangle  ,
$$
whenever $\phi\in\scr D$, and $f$ belongs to $\scr D'$.  In view of the skewness,
we define 
$$\frac{\partial f}{\partial\bar z}
=
-
\left(\frac{\partial}{\partial\bar z}\right)^*f, \ \forall f\in\scr D',
$$
so that the operator $\frac{\partial}{\partial\bar z}$ on $\scr D'$
is the weak-star continuous extension of   
$\frac{\partial}{\partial\bar z}$ on $\scr D$, when we regard $\scr D$
as a subset of $\scr D'$, under the identification of each $f\in L^1_{\textup{loc}}(dxdy)$ with the distribution {\em represented by $f$}, given by
$$
\langle\phi,f\rangle =
\int_{\mathbb{C}}\phi f dxdy,\ \forall\phi\in\scr D.$$

\medskip
Complex Radon measures (Borel-regular complex-valued measures on $\mathbb{C}$,
having finite total variation on each compact subset of $\mathbb{C}$) also
represent distributions. The measure $\mu$ acts continuously on the space
$C^0_{\textup{cs}}$ of continuous complex-valued
functions on $\mathbb{C}$ having compact support, and equipped with the usual inductive limit topology, via
$$ \langle\phi,\mu\rangle =
\int_{\mathbb{C}}\phi f d\mu,\ \forall\phi\in C^0_{\textup{cs}}$$
and hence restricts to a continuous linear functional on $\scr D$.
If we identify $f\in L^1_{\textup{loc}}$ with the measure $fdxdy$, then 
this generalises the previous remark.
A measure is uniquely-determined by the corresponding distribution, because 
$\scr D$ is dense in $C^0_{\textup{cs}}$, so we may identify the measure
and distribution, without fear of confusion.

\medskip
The $\bar\partial$ operator is linear and translation-invariant. It is also
{\em elliptic}: this means that it is almost invertible; more precisely it has
finite-dimensional kernel and cokernel, when restricted to a suitable space.
When restricted to $\scr D'$, it has a very big kernel, the space of all entire
functions. This statement is a case of Weyl's Lemma: a distributional 
solution $u\in\scr D'(U)$
(where $\scr D(U)$ is the space of $C^\infty$ complex-valued functions on $U$)
of $ \bar\partial u=0$ on an open set $U\subset\mathbb{C}$ is representable
by a holomorphic function on $U$. This is a good thing, because it gives
complex analysts a nontrivial field of study. But when restricted to
$\scr E'$, the dual of the space of $C^\infty$ functions on $\mathbb{C}$
with compact support,
$\bar\partial$ is injective. The fundamental solution is
the locally-integrable function $-1/\pi z$, i.e.
$$ \frac{\partial}{\partial\bar z}\left(-\frac1{\pi z}\right) = \delta_0,$$
the point mass at $0$. For $\phi\in\scr D$, we have
$$ \frac{\partial}{\partial\bar z}\hat\phi = \phi,$$
where 
$$\hat\phi = 
\left(-\frac1{\pi z}\right)*\phi,$$
the {\em Cauchy transform} of $\phi$. 
We extend the transform to a map $\scr E'\to\scr D'$ by setting
$$\langle\phi,\hat f\rangle =-\langle\hat\phi,f\rangle,
\ \forall\phi\in\scr D,\, \forall f\in\scr E'.$$
We have
$$ \frac{\partial}{\partial\bar z}\hat f = f,\ \forall f\in\scr E'.$$
In particular, for each $f\in\scr E'$, the distribution $\hat f$
is (represented by) a holomorphic function off $\textup{spt}f$,
(the support of $f$.
Because of all this, the Cauchy transform is intimately
connected with analytic function theory, and one can use it
to establish many interesting results.

\section{Some Holomorphic Approximation Theorems}
For $X\subset\mathbb{C}^n$, let $\scr O(X)$ denote the space
of functions holomorphic near $X$. For compact Hausdorff
$X$, let $C^0(X)$ denote the Banach space of all continuous,
complex-valued functions on $X$, with the sup norm. 
\begin{theorem}[Hartogs-Rosenthal, 1931] 
Suppose $X\subset\mathbb{C}$ is compact and has area zero. Then 
$\scr O(X)$ is dense in $C^0(X)$.
\end{theorem}
\begin{proof} By the Separation Theorem for Banach spaces,
it suffices to show that 
$$ L\in C^0(X)^*\cap\scr O(X)^{\perp} \Rightarrow L=0.$$
By the Riesz Representation Theorem, the dual $C^0(X)^*=M(X)$,
the space of (complex, Radon) measures supported on $X$.

Fix $\mu\in M(X)$ with $\mu\perp\scr O(X)$, i.e.
$\int f\,d\mu=0$ whenever $f\in\scr O(X)$.

Regarding $\mu$ as a distribution on $\mathbb{C}$, we find that
$\hat\mu$ is the locally-integrable function given by
$$ \hat\mu(\zeta) = \frac1{\pi}\int\frac{d\mu(z)}{z-\zeta},
\ \forall\zeta\in\mathbb{C}.$$
Since the function $z\mapsto1/(z-\zeta)$ belongs to $\scr O(X)$
for $\zeta\not\in X$, we have $\hat\mu=0$ $dxdy$-a.e., hence $\hat\mu=0$
as a distribution, hence 
$$\mu = \frac{\partial\hat\mu}{\partial\bar z}=0.$$
\end{proof}

\begin{corollary}[A. Browder] Suppose $X\subset\mathbb{C}^n$ is compact,
and each coordinate 
projection of $X$ has area zero. Then $\scr O(X)$ is dense in $C^0(X)$.
\end{corollary}
\begin{proof}
Denote $z=(z_1,\ldots,z_n)$, and $z_j=x_j+iy_j$. 
Fix $j\in\{1,\ldots,n\}$.
Let $\pi_j:z\mapsto z_j$. Then $\pi_j(X)$
has area zero, so by Hartogs-Rosenthal 
$z\mapsto x_j$ and $z\mapsto y_j$
are uniform limits of functions (depending only on $\pi_j(z)$)
that are holomorphic on a neighbourhood of $\pi_j^{-1}(\pi_j(X))$,
and hence belong to $\scr O(X)$. 
Thus the uniform closure $A$ on $X$ of the algebra $\scr O(X)$
contains all the coordinate functions $x_j$ and $y_j$,
so by La Valle\'e Poussin's extension of Weierstrass'
Polynomial Approximation Theorem (a special case of the
Stone-Weierstrass Theorem), we conclude that
$A=C^0(X)$.
\end{proof}

\medskip
Consider the function spaces, for $0<\alpha<1$ and compact
$X\subset\mathbb{C}^n$:
$$\textup{Lip}(\alpha,X) = \left\{
f\in C^0(X): 
\sup_{z\not=w}\frac{|f(z)-f(w)|}{|z-w|^\alpha}<+\infty\right\},$$
and
$$\textup{lip}(\alpha,X) = \left\{
f\in\textup{Lip}(\alpha,X): 
\sup_{0<|z-w|<\delta}\frac{|f(z)-f(w)|}{|z-w|^\alpha}\to0\textup{ as }\delta\downarrow0\right\}.$$
With a suitable norm, $\textup{Lip}(\alpha,X)$ becomes a Banach algebra, and 
the subspace $\textup{lip}(\alpha,X)$ is a closed subalgebra, equal to the 
closure of $\scr D$ in $\textup{Lip}(\alpha,X)$.
The elements of the dual $\textup{lip}(\alpha,X)^*$ may be represented
in a manner somewhat similar to the Riesz representation, as follows.

Fix any $a_0\in X$.

Given $L\in\textup{lip}(\alpha,X)^*$, there exist $\lambda\in\mathbb{C}$ and a
measure $\mu$ on the product $X\times X$ having no mass on the diagonal, such that
$$
Lf = \lambda f(a_0) +
\int_{X\times X}\frac{f(z)-f(w)}{|z-w|^\alpha}
d\mu(z,w),$$
whenever $f\in\textup{lip}(\alpha,X)$.  If $L1=0$, then $\lambda=0$.
For such $\lambda$, this permits us to represent
$\hat L$ by integration against an $L^1_{\textup{loc}}$ function:
$$ \hat L(\zeta) = \frac1{\pi}\int\frac{(w-z)\,d\mu(z,w)}{%
(\zeta-z)(\zeta-w)|z-w|^\alpha}.$$

Using this, and essentially the same proof as given for
the Hartogs-Rosenthal Theorem, one obtains \cite[p. 387]{OF1}:

\begin{theorem} If $X\subset\mathbb{C}$ is compact with area
zero, then $\scr O(X)$ is dense in $\textup{lip}(\alpha,X)$
for $0<\alpha<1$.
\end{theorem}
\begin{corollary}
If $X\subset\mathbb{C}^n$ has all its coordinate projections of
area zero, then $\scr O(X)$ is dense in $\textup{lip}(\alpha,X)$.
\end{corollary}

\section{Higher-dimensional Cauchy Transforms}
The utility of the Cauchy transform in one dimension prompted
people to seek a similar tool for problems of several complex
variables. Here one {\em must} use forms. The kernel
$-1/\pi z$ must be replaced by a $(2n-1)$-form of type $(n,n-1)$:
$$
\Omega = \sum_{j=1}^n K_j(\zeta,z)\,
d\bar\zeta_1\wedge\cdots d\bar\zeta_{j-1}\wedge d\bar\zeta_{j+1}\wedge
\cdots\wedge d\bar\zeta_n\wedge d\zeta_1\wedge\cdots\wedge
d\zeta_n,$$
such that
\begin{equation}\label{eqn-*}
 \phi(z) = \int\Omega(\zeta,z)\wedge\bar\partial\phi(\zeta)
\end{equation}
holds for test functions $\phi$. A form $\Omega$
that does this is called a {\em Cauchy-Leray-Fantappi\'e form}.
There are many such forms, and depending on the end in view, one
prefers one or another.  There are also more complex forms,
involving boundary terms (analogous to Pompeiu's formula), useful
for specific purposes.

In joint work with David Walsh, and the late 
Ken Preskenis, we obtained the following \cite{OPW1984}:
\begin{theorem} Let $X\subset\mathbb{C}^n$ be compact and holomorphically-convex.Let $E\subset X$ be closed, and suppose that each point $a\in X\sim E$ has a neighbourhood $N\subset\mathbb{C}^n$ such that
$X\cap N$ is a subset of a $C^1$ submanifold without complex tangents. 
Then
$$ \textup{clos}_{C^0(X)}\scr O(X) = C^0(X)\cap
\textup{clos}_{C^0(E)}\scr O(X),$$
and
$$ \textup{clos}_{\textup{Lip}(\alpha,X)}\scr O(X) = \textup{lip}(\alpha,X)\cap
\textup{clos}_{\textup{Lip}(\alpha,E)}\scr O(X),\ \textup{ for }0<\alpha<1.$$
\end{theorem}
In other words, approximation problems on $X$ reduce to approximation
problems on the singular set $E\subset X$. 

This generalised and extended to Lip$(\alpha)$
earlier work of Range and Siu ($E=\emptyset$),
Weinstock ($X$ polynomially-convex), and ourselves \cite{OPW1983}
(See below). 

The proof comes down to showing that if a distribution 
$L$ that acts continuously on $\textup{lip}(\alpha,X)^*$ 
annihilates $\scr O(X)$, then $L$
is supported on $E$. To do this, one constructs a kernel $\Omega(\zeta,z)$
such that Equation (\ref{eqn-*}) holds for $z$ on a neighbourhood 
$U$ of $X$ and $\phi\in\scr D(U)$, and a second kernel
$\tilde\Omega(\zeta,z)$ such that $\tilde\Omega(\zeta,z)=\Omega(\zeta,z)$
for $z\in X$ and $\zeta\in U$, and $\tilde\Omega$ has
coefficients $\tilde K_j(\zeta,z)$ that are holomorphic
in $z\in U$ for each $\zeta\in U\sim X$, and have another
technical property. This construction is based on work of
Berndtsson, building on the special Bochner-Martinelli kernel.
Then, representing $L$ as before by a measure $\mu$ on
$X\times X$ with no mass on the diagonal, we can represent
$$
\begin{array}{rcl}
&\langle\phi,& L\rangle\\
&=&
\displaystyle\int_{X\times X}\displaystyle\frac1{|z-w|^\alpha}
\int_U\left\{\Omega(\zeta,z)-\Omega(\zeta,w)\right\}
\wedge \bar\partial\phi(\zeta)\,d\mu(z,w)\\
&=&
\displaystyle\int_U
\displaystyle\int_{X\times X}\displaystyle\frac1{|z-w|^\alpha}
\left\{\Omega(\zeta,z)-\Omega(\zeta,w)\right\}
d\mu(z,w)
\wedge \bar\partial\phi(\zeta),\\
\end{array}
$$
by Fubini's Theorem. (There are substantial technical estimates
involved in justifying this.)

It remains to show that
$$  \int_{X\times X}\frac1{|z-w|^\alpha}
\left\{\Omega(\zeta,z)-\Omega(\zeta,w)\right\}
d\mu(z,w)=0
$$
for almost all $\zeta\in U$. The fact that 
$\Omega(\zeta,z)=
\tilde\Omega(\zeta,z)$
for $z\in X$ and that the latter is holomorphic in  
$z$, and the technical properties (the most important of which is
an \lq\lq omitted sector property") allow us to approximate
each coefficient in the integral by elements of $\scr O(X)$,
and gives the desired result. For the details, see
\cite{OPW1984}.

\begin{corollary}[Range-Siu] If $E=\emptyset$, then $\scr O(X)$
is dense in $C^0(X)$.
\end{corollary}
\begin{corollary}
Let $F\subset Y$, where $Y$ is a compact subset of $\mathbb{C}^n$
and $F$ is a closed subset of $Y$. Let $f$ be a $\mathbb{C}^n$-valued
function defined on a neighbourhood of $Y$, let $X=f(Y)$
and $E=f(F)$. Suppose that $X$ is polynomially-convex, and
the matrix $f_{\bar z}$ (with columns $\frac{\partial f}{\partial\bar z_j}$)
is invertible on $Y\sim F$. Then
\begin{equation}\label{eqn-1}
\textup{clos}_{C^0(X)}\mathbb{C}[z,w] =
C^0(X)\cap 
\textup{clos}_{C^0(E)}\mathbb{C}[z,w],
\end{equation}
and
\begin{equation}\label{eqn-2}
\textup{clos}_{\textup{Lip}(\alpha,X)}\mathbb{C}[z,w] =
\textup{lip}(\alpha,X)\cap 
\textup{clos}_{\textup{Lip}(\alpha,E)}\mathbb{C}[z,w].
\end{equation} 
\end{corollary}
Equation \ref{eqn-1} is due to Weinstock.
\begin{corollary}
Suppose $\rho$ is a $C^2$ strictly plurisubharmonic function on
a neighbourhood of $\textup{bdy}X$, where $X$ is a compact subset of 
$\mathbb{C}^n$, with interior $D$, and that
$\textup{bdy}X=\{z: \rho(z)=0\}$, 
$\{z:\rho(z)<0\}\subset D$, and $E=\textup{clos}D$.
Then Equations \ref{eqn-1} and \ref{eqn-2} hold.
\end{corollary}
In this case, Equation \ref{eqn-1} is due to Henkin and Leiterer.

\section{Extending Smooth Functions}
According to one view, Geometry is an aspect of Group Theory.
But more accurately, Geometry {\em is}
Ring Theory. To be absolutely precise, Geometry
is Topological Ring Theory.

Let $M$ be a $C^k$ manifold, and $X\subset M$. Given
$f:X\to \mathbb{R}$, when does there exist a $C^k$ function
$\tilde f:M\to\mathbb{R}$ such that the restriction $\tilde f|X=f$?
This problem arises in many applications,
and has been studied since the 1930's, with important work of
Whitney and Glaeser. Richard Watson and I studied it in the
early 1990's drawing on some ideas of mine that go back to the 
1970's. 

We deal now with real-valued functions, real vector spaces
and algebras, and $i$ is just an index, or multi-index,
and not $\sqrt{-1}$
any more.

Let $C^k(M)$ now denote the algebra (under pointwise operations)
of $C^k$ real-valued functions on $M$. This is a Frechet algebra
(a complete metric algebra) with the natural topology. For
$S\subset C^k(M)$, let
$$ S^\perp = \{L\in C^k(M)^*: Lf=0,\forall f\in S\},$$
where $C^k(M)^*$ denotes the space of continuous linear
functionals $L:C^k(M)\to\mathbb{R}$. Note that $C^k(M)^*$
is a module over $C^k(M)$. For $X\subset M$, let
$$ X_\perp =\{f\in C^k(M): f(a)=0,\forall a\in X\}.$$
Let $a_\perp = \{a\}_\perp$, when $a\in X$. Each $X_\perp$
is an ideal in $C^k(M)$. The ideal $(a_\perp)^{k+1}$ is generated by
products of $k+1$ elements of $a_\perp$. Its annihilator
$((a_\perp)^{k+1})^\perp$ consists of the so-called
\lq\lq $k$-th order point differential operators". In local coordinates
$(x_1,\ldots,x_d)$, each $\partial\in
((a_\perp)^{k+1})^\perp$ takes the form
$$ \partial f = \sum_{|i|\le k}\alpha_i
\frac{\partial f}{\partial x_i}(a),\ \forall f\in C^k(M),$$
where $i=(i_1,\ldots,i_d)\in\mathbb{Z}_+^d$ denotes a 
multi=index,
$|i|=\sum_j i_j$, and
$\alpha_i\in\mathbb{R}$ are constants depending on $\partial$,
but not on $f$.

The {\em $k$-th order tangent space to $M$ at $a$} is
defined as
$$ \textup{Tan}^k(M,a) = C^k(M)^*\cap 
((a_\perp)^{k+1})^\perp,$$
and the 
{\em $k$-th order tangent space to $X$ at $a$} is
defined as
$$ \textup{Tan}^k(M,X,a) =  
 \textup{Tan}^k(M,a)\cap(X_\perp)^\perp,$$
the set of $k$-th order point differential operators $\partial$ at $a$
such that $\partial f$ depends only on the values of $f$ on $X$.  The disjoint
unions
$$\begin{array}{rcl}
T^k(M) &=& \dot\bigcup_{a\in M} \textup{Tan}^k(M,a),\\
\textup{and}&&\\
T^k(M,X) &=& \dot\bigcup_{a\in M} \textup{Tan}^k(M,X,a)
\end{array}
$$
are called the {\em $k$-th order tangent bundle of 
$M$}, and the {\em $k$-th order tangent sheaf of $X$},
respectively.  The stalks $\textup{Tan}^k(M,X,a)$ have the structure
of finite-dimensional modules over a finite-dimensional real algebra, and 
provide numerical $C^k$ invariants for the pair $(M,X)$, since the 
tangent construction behaves functorially. A $C^k$ function
$F:M\to M'$ between $C^k$ manifolds induces
an algebra homomorphism
$$ F^\#:
\left\{
\begin{array}{rcl}
C^k(M') &\to& C^k(M)\\
g&\mapsto& g\circ F
\end{array}
\right.
$$
and a $C^k$-module homomorphism
$$F_\# = (F^\#)^*: C^k(M)^*\to C^k(M')^*.$$
If $F$ maps $X$ into $X'$, then $F_\#$ maps 
the stalk $\textup{Tan}^k(M,X,a)$ to the stalk 
$\textup{Tan}^k(M',X',f(a))$, and so
induces a map $F_*:T^k(M,X)\to T^k(M',X')$.
We established the following \cite{OW}:
\begin{theorem} Let $X$ be a closed subset of a 
$C^k$ manifold $M$, $f:X\to\mathbb{R}$ be continuous, and
$$
\pi:\left\{
\begin{array}{rcl}
M\times\mathbb{R} &\to& M\\
(x,y)&\mapsto& x
\end{array}
\right.
$$
be the projection. Then $f$ has a $C^k$ extension to $M$
if and only if the map
$$ \pi_*: T^k(M\times \mathbb{R},f) \to T^k(M,X)$$
is bijective.
\end{theorem}
This result, and the $k$-th order tangent concept, are not
particularly difficult, but are completely fundamental for the
extension problem.  They reduce the extension problem
to the problem of deciding whether or not two integral
dimensions (of $\textup{Tan}^k(M\times\mathbb{R},f,(a,f(a))$
and $\textup{Tan}^k(M,X,a)$) agree at each point $a\in X$.
We were gratified by the favourable reception of
our paper, which included a congratulatory
letter from Malgrange. 
It remained a problem to come up with a constructive procedure
for deciding the question. Whitney himself dealt with this in dimension
one, for all $k$. We provided a way to do it
in 1993, in case $k=1$, for all dimensions. In recent years, 
C. Fefferman and co-workers have gone as far as can be
done in providing a constructive procedure for general $k$. 
See the website \cite{F}
where this monumental corpus may be downloaded.
Their
work employs, {\em inter alia} the $k$-th order sheaf we introduced
and our result. Fefferman was unaware of our 
work, having taken the concept and result from the 2003 Inventiones
paper of Bierstone, Milman and Pawluckii \cite{BMP}. 
I supplied a copy of our
paper to Pawlucki in December 1997, at his request. These authors
included our paper among their references, but did not attribute
the concept to us. They referred to our paper only in order
to make a gratuitously dismissive remark about it. I am at a loss
to understand this behaviour. They called the tangent sheaf
the {\em Zariski paratangent bundle} (although it is not in general
a bundle), and subsequently a 
number of authors have referred to it as {\em the Zariski paratangent
bundle  of Bierstone, Milman and Pawlucki}.  The bundle $T^k(M)$
was originally introduced by Pohl and Feldman, but the sheaf
$T^k(M,X)$ appeared first in our paper.  

I must also mention the work of Declan O'Keeffe \cite{DOK}, 
who proved the analogous
result for $C^{k+\alpha}$ extensions ($k\in\mathbb{N}$, $0<\alpha<1$),
and who also used $T^k$ to study algebraic curve
singularities in $\mathbb{C}^2$.

\section{Approximating $C^\infty$ Functions}
In conclusion, here is a brief summary of the joint work with the late Graham Allan, Grayson Kakiko and Richard, on Segal's problem. The problem called for a characterization
of the closed subalgebras of the algebra $C^\infty(M,\mathbb{R})$, for a smooth manifold
$M$. The problem is local, so that we may take $M=\mathbb{R}^d$.  There is not
much loss in generality in considering subalgebras that are 
topologically-finitely-generated, i.e. those of the form
$$ A(\Psi) = \textup{clos}_{C^\infty(\mathbb{R}^d)}\mathbb{R}[\Psi]
= 
\textup{clos}_{C^\infty(\mathbb{R}^d)}\{g\circ\Psi:g\in C^\infty(\mathbb{R}^r,
\mathbb{R}),$$
where $\Psi=(\psi_1,\ldots,\psi_r)\in C^\infty(\mathbb{R}^d$. In 1950, 
Nachbin conjectured a solution, analogous to Whitney's Spectral Theorem for
closed ideals. Let $\mathbb{R}[[x_1,\ldots,x_n]]$ denote the algebra
of formal power series in $n$ indeterminates, let
$$ T'_a: C^\infty(\mathbb{R}^d,\mathbb{R}^r) \to
\mathbb{R}[[x_1,\ldots,x_d]]^r$$
denote the truncated Taylor series map, for each $a\in\mathbb{R}^d$,
and let $T_af=f(a)+T'_af$ be the full Taylor series. (Note that
$x_j$ would have to be replaced by $(x_j-a_j)$ in these series
in applications of Taylor's Theorem.)

In case $\Psi\in C^\infty(\mathbb{R}^d,\mathbb{R}^r)$ is injective,
Nachbin's conjecture comes down to
\begin{equation}\label{eqn-**}
A(\Psi) = \bigcap_{a\in\textup{crit}\Psi}
T_a^{-1}\mathbb{R}[[T'_a\Psi]],
\end{equation}
which may be stated in loose terms as: $f\in A(\Psi)$ if and only if
$f$ has the \lq\lq right kind" of Taylor series at each point.

\begin{theorem}[Tougeron 1971]
Suppose that for each compact $K\subset\mathbb{R}^d$ there exists
$\alpha>0$ and $\beta>0$ such that
$$|\Psi(x)-\Psi(y)| \ge \alpha|x-y|^\beta,\ \forall x,y\in K.$$
Then Equation (\ref{eqn-**}) holds.
\end{theorem}

\begin{corollary}If $\Psi$ is injective and real-analytic, then
Equation (\ref{eqn-**}) holds.
\end{corollary}

Our main result was this \cite{AKOW1998a}:
\begin{theorem} If $\Psi$ is injective and $d=1$, then Equation 
(\ref{eqn-**}) holds.
\end{theorem}

The proof involves some hard analysis, to oversome the problems around
the accumulation points of the critical set. 

We say that $\Psi$ is {\em flat at $a$} if $T'_a\Psi=0$.
We also proved the following useful result \cite{AKOW1998b}
\begin{theorem} If $\Psi$ is injective and $f\in C^\infty(\mathbb{R}^d)$
is flat at each critical point of $\Psi$, then $f\in A(\Psi)$.
\end{theorem}

\begin{corollary} If $\Psi$ is injective and flat on $\textup{crit}\Psi$,
then Equation (\ref{eqn-**}) holds.
\end{corollary}
\begin{corollary} If $\Psi$ is injective and $\textup{crit}\Psi$ is discrete,
then Equation (\ref{eqn-**}) holds.
\end{corollary}
In further work \cite{AKOW2001}, we studied $A(\Psi)$ as a Frechet algebra, for injective $\Psi$. We established that it is always regular, and that membership
of $A(\Psi)$ is a local property.

\section{Notes}
It is convenient to take this opportunity to correct 
the definition of proxy distance given on p.49 of our paper
\cite{AKOW1998b} in the proceedings of the meeting at
Blaubeuren. This should read:

{\bf Definition} Let $E\subset\mathbb{R}^d$ be closed and
$\kappa_m\ge1$ ($m=0,1,2,\ldots$).
A function $d_E:\mathbb{R}^d\to [0,+\infty)$ that
is $C^\infty$ on $\mathbb{R}^di\sim E$ is
called a $\{\kappa_m\}$ {\em proxy distance for $E$}
if 
$$ \frac1{\kappa_0}d_E(x)\le\textup{dist}(x,E)\le\kappa_0d_E(x),\ \forall
x\in\mathbb{R}^d$$
and
$$
|D^md_E(x)| \le \kappa_m\cdot\textup{dist}(x,E)^{1-m},\ \forall x\in\mathbb{R}^d,
\forall m\ge1.$$

No other change to the paper is needed, and it all remains true.

(The reason for the change is that $d_E$ {\em cannot} be $C^\infty$
on the whole of $\mathbb{R}^d$ when $\emptyset\not= E\not=\mathbb{R}^d$;
in the subsequent application $d_E$ is always composed with functions 
$\phi$ that vanish near $0$, so $\phi\circ d_E$ is
$C^\infty$ on $\mathbb{R}^d$.)

\medskip
I would also like to point out to
readers of \cite{OPW1984}
that the interesting case of Example 3.2 is 
when
$X\not=\textup{int}\,\textup{clos}\,X$. It is even
interesting when $X$ has no interior.

\end{document}